\documentclass[11pt,reqno]{amsart}

\usepackage{amssymb, amsmath}
\usepackage{hyperref}
\usepackage[nohug]{diagrams}
\diagramstyle[labelstyle=\scriptstyle]

\newtheorem{proposition}{Proposition}[section]
\newtheorem{lemma}[proposition]{Lemma}
\newtheorem{corollary}[proposition]{Corollary}
\newtheorem{theorem}[proposition]{Theorem}

\theoremstyle{definition}

\newtheorem{example}[proposition]{Example}

\newtheorem{remark}[proposition]{Remark}

\newcommand{\thlabel}[1]{\label{th:#1}}
\newcommand{\thref}[1]{Theorem~\ref{th:#1}}
\newcommand{\selabel}[1]{\label{se:#1}}
\newcommand{\seref}[1]{Section~\ref{se:#1}}
\newcommand{\lelabel}[1]{\label{le:#1}}
\newcommand{\leref}[1]{Lemma~\ref{le:#1}}
\newcommand{\prlabel}[1]{\label{pr:#1}}
\newcommand{\prref}[1]{Proposition~\ref{pr:#1}}
\newcommand{\colabel}[1]{\label{co:#1}}
\newcommand{\coref}[1]{Corollary~\ref{co:#1}}
\newcommand{\relabel}[1]{\label{re:#1}}
\newcommand{\reref}[1]{Remark~\ref{re:#1}}
\newcommand{\exlabel}[1]{\label{ex:#1}}

\newcommand{\eqlabel}[1]{\label{eq:#1}}
\newcommand{\equref}[1]{(\ref{eq:#1})}

\newcommand{\mc}{\mathcal}
\newcommand{\mb}{\mathbb}

\newcommand{\ve}{\varepsilon}

\newcommand{\ol}{\overline}
\newcommand{\la}{\langle}  \newcommand{\ra}{\rangle}

\newcommand{\coalg}{{\rm CoAlg}}

\newcommand{\halg}{{\rm HopfAlg}}

\begin{document}

\title{Subcoalgebras and endomorphisms of free Hopf algebras}
\author{Alexandru Chirv\u asitu}

\address{
University of California, Berkeley, 970 Evans Hall \#3480, Berkeley CA, 94720-3840, USA
}
\email{chirvasitua@gmail.com}


\begin{abstract}

For a matrix coalgebra $C$ over some field, we determine all small subcoalgebras of the free Hopf algebra on $C$, the free Hopf algebra with a bijective antipode on $C$, and the free Hopf algebra with antipode $S$ satisfying $S^{2d}={\rm id}$ on $C$ for some fixed $d$. We use this information to find the endomorphisms of these free Hopf algebras, and to determine the centers of the categories of Hopf algebras, Hopf algebras with bijective antipode, and Hopf algebras with antipode of order dividing $2d$. 

\end{abstract}

\maketitle

\section*{Introduction}\selabel{0}

The free Hopf algebra $H(C)$ on a coalgebra $C$ (over some base field $k$) was introduced by Takeuchi in \cite{Ta}, and this construction was used to give the first example of a Hopf algebra with non-bijective antipode: if $n>1$, Takeuchi shows that the antipode of $H(M_n(k)^*)$ is not bijective. Later, Nichols (\cite{Ni}) constructed bases for $H(M_n(k)^*)$, and showed that for $n>1$ the antipode is injective. 

In \cite{Sc}, Schauenburg constructed the free Hopf algebra with bijective antipode on a Hopf algebra. This can be combined with Takeuchi's construction to yield a left adjoint for the forgetful functor from the category of Hopf algebras with bijective antipode to that of coalgebras (\cite[Lemma 3.1]{Sc}). In the same paper, a basis for the free Hopf algebra with bijective antipode (which we will denote by $H_\infty(C)$) on a matrix coalgebra $C$ was constructed, by methods analogous to those used by Nichols in \cite{Ni}, and used to give the first examples of a Hopf algebra with surjective, non-injective antipode.

In this paper we study the subcoalgebras and endomorphisms of these objects, and also of the free Hopf algebra $H_d(M_n(k)^*)$ with antipode whose order divides $2d$ ($n>1$). This seems not to have been done in too much detail in the literature. Further motivation comes from the desire to study the so-called {\it centers} of the categories appearing in the above discussion, i.e. their monoids of self-natural transformations of the identity functor: the category $\halg$ of Hopf algebras, $\halg_\infty$ of Hopf algebras with bijective antipode, and $\halg_d$ of Hopf algebras with antipode of a given order $2d$ (for some positive integer $d$). 

For a justification for the term ``center'', notice that if our category is a monoid (i.e. a one-object category), then its center is precisely the center of the monoid. As another example, notice that there is an obvious isomorphism between the center of the category $_A\mc M$ of left modules over a ring $A$, and the center of $A$.  

One would expect, for example, that the center of $\halg$ is the free monoid on one element generated by the square of the antipode, together with a ``multiplicative $0$'', the (natural transformation induced by the) trivial endomorphism. It is obvious from this statement that antipodes cannot all be bijective, so it can be regarded as a natural generalization of Takeuchi's results in \cite{Ta}. We prove this result and the analogous ones for $\halg_\infty$ and $\halg_d$ below, using the fact that, as will become apparent, the free Hopf algebras mentioned above on an $n\times n$ matrix coalgebra ($n>1$) have, in a certain sense, ``no more endomorphisms than expected'' (in most cases).

The paper is organized as follows:

In \seref{1} we introduce the notations and conventions to be used throughout, and recall a few facts on free Hopf algebras. 

In \seref{2} the main results are proven. We show that the free Hopf algebra $H(C)$, the free Hopf algebra with bijective antipode $H_\infty(C)$, and the free Hopf algebra $H_d(C)$ with antipode of order $2d\ge 4$ on an $n\times n$ matrix coalgebra $C$ ($n>1$) contains no subcoalgebras of dimension $\le n^2$, other than the $1$-dimensional coalgebra $k$ and the iterates of $C$ by the antipode. $d=1$ is a little trickier, so in this case, we prove our result only for $n>2$ or in characteristic zero. We then use this to prove that the Hopf algebra endomorphisms of $H(C),\ H_\infty(C)$ and $H_d(C)$ are precisely those we would expect, in most cases (i.e. apart from $d=1,n=2$, positive characteristic): the compositions of those induced by the (anti)endomorphisms of $C$ with the powers of the antipode, and also the trivial endomorphism (unit composed with counit); note that an {\it anti}endomorphism of $C$ is required if we are to compose with an odd power of the antipode. 

In \seref{3}, the results outlined above are used to determine the centers of the categories $\halg,\ \halg_\infty$, and $\halg_d$. Again, the final result is exactly as expected: in all three cases we have the (natural transformations induced by) the trivial endomorphism, and the even powers $S^{2t}$ of the antipode $S$. Of course, in the three cases in question, $t$ ranges through the appropriate sets: the non-negative integers, the integers, and $\mb Z/2d$ respectively. 

Finally, in \seref{4} we look at the exceptions mentioned before: $d=1,n=2$, positive characteristic. It is shown that indeed, there are counterexamples to the results in \seref{2} in characteristic $2$ or $3$. 

\section{Preliminaries}\selabel{1}

We work over some fixed field $k$. Algebras, coalgebras, Hopf algebras, etc. are over $k$, and all (co)algebras are (co)unital and (co)associative. We assume familiarity with basic Hopf algebra theory, as in \cite{Sw,A} or \cite{Mo}, for instance. We denote the categories of Hopf algebras, Hopf algebras with bijective antipode, and Hopf algebras with antipode $S$ such that $S^{2d}={\rm id}$ ($d\ge 1$) by $\halg,\ \halg_\infty$ and $\halg_d$ respectively. We reserve the usual notations for other categories that might appear ($\coalg$ is the category of $k$-coalgebras, for example). The usual symbols are used for the structure maps of our coalgebras, bialgebras, etc.: $\Delta,\ve, S$ denote the comultiplication, counit and respectively the antipode of an appropriate object. We might use the name of the object as a subscript for the structure map: $S_H$ is the antipode of the Hopf algebra $H$, for example. 

Recall (\cite{Ta}) that the forgetful functor $\halg\to\coalg$ has a left adjoint; there is a free Hopf algebra $H(C)$ on any coalgebra $C$, with the usual universal property. Similarly (\cite{Sc}), there is a free Hopf algebra $H_\infty(C)$ on any coalgebra $C$. Using the exact same techniques as in those papers, or, alternatively, just factoring $H_\infty(C)$ through the appropriate ideal, we have the following result:

\begin{proposition}\prlabel{H_d}

For every positive integer $d$, the forgetful functor $\halg_d\to\coalg$ has a left adjoint. 
 
\end{proposition}

We denote this left adjoint by $H_d(-)$. The proof is entirely routine, and is left to the reader; one simply factors $H_\infty(C)$ (or even $H(C)$) through the appropriate ideal to get $H_d(C)$. When we wish to state a result in a unified manner for $H,H_\infty$ and $H_d$ all at once, we use the notation $\tilde H(C)$ to stand for either one of them. 

In most cases for us (but not always), $C$ will be an $n\times n$ matrix coalgebra for some $n>1$. This is the dual $M_n(k)^*$ of the matrix algebra $M_n(k)$; it has a basis $(t_{ij})_{i,j=1}^n$, with the coalgebra structure given by
\begin{equation}\eqlabel{matco}
\Delta(t_{ij})=\sum_k t_{ik}\otimes t_{kj},\ \ve(t_{ij})=\delta_{ij},
\end{equation}
where $\delta_{ij}$ is the Kronecker delta, as usual. 

An important tool for us will be the $k$-basis constructed by Nichols for $H(M_n(k)^*)$ (\cite{Ni}), and the analogous bases for $H_\infty$ (\cite{Sc}) and $H_d$. To my knowledge, the $H_d$ case has not appeared in the literature, but the calculations are mostly parallel to those used for $H$ and $H_\infty$, and we do not repeat those here. The only problem when trying to adapt the proof in \cite{Ni} to $H_d$ arises when $d=1$. We address this briefly below, again, omitting the verifications.   

Recall the notation $\tilde H$ introduced above. Below, it is understood that $r$ ranges through the non-negative integers if $\tilde H=H$, through $\mb Z$ if $\tilde H=H_\infty$, and through $\mb Z/2d$ if $\tilde H=H_d$. Consider the set $\mc X=\{x_{ij}^r,\ |\ i,j=\ol{1,n},r\}$. We now work inside the free algebra $k\la \mc X\ra$ on $\mc X$, and seek to write $\tilde H(M_n(k)^*)$ as a quotient of $k\la\mc X\ra$, using Bergman's diamond lemma (\cite{Be}). The images of $x^0_{ij}$ in $\tilde H(M_n(k)^*)$ are supposed to be matrix generators for $M_n(k)^*$, and we want the antipode to act by sending $x^r_{ij}$ to $x^{r+1}_{ji}$. The diamond lemma comes in when trying to factor out the relations imposed by the condition that this map be an antipode. In \cite{Ni, Sc}, this is done as follows (for $\tilde H=H$ and $H_\infty$, respectively):  

Following \cite{Sc}, we say that a monomial $w$ in the free monoid $\la \mc X\ra$ on $\mc X$ is less than another monomial $w'$ if either $w$ is shorter than $w'$, or if they have the same length, the same sequence of $r$-indices, and the sequence of $i,j$ indices of $w$ is lexicographically less than that of $w'$. Now consider the reductions
\begin{equation}\eqlabel{red1}
x_{in}^r x_{jn}^{r+1}     \to    \delta_{ij}-\sum_{a<n}x_{ia}^r x_{ja}^{r+1}
\end{equation}
\begin{equation}\eqlabel{red2}
x_{ni}^{r+1} x_{nj}^r     \to    \delta_{ij}-\sum_{a<n}x_{ai}^{r+1} x_{aj}^r    
\end{equation}
\begin{equation}\eqlabel{red3}
x_{in}^r x_{jn-1}^{r+1} x_{kn-1}^{r+2}    \to    \delta_{jk} x_{in}^r - \delta_{ij}x_{kn}^{r+2}+\sum_{a<n}x_{ia}^r x_{ja}^{r+1} x_{kn}^{r+2}																														-\sum_{a<n-1} x_{in}^r x_{ja}^{r+1} x_{ka}^{r+2} 																						
\end{equation}
\begin{equation}\eqlabel{red4}
x_{ni}^{r+2} x_{n-1j}^{r+1} x_{n-1k}^r    \to    \delta_{jk} x_{ni}^{r+2} - \delta_{ij}x_{nk}^r+\sum_{a<n}x_{ai}^{r+2} x_{aj}^{r+1} x_{nk}^r
																								-\sum_{a<n-1} x_{ni}^{r+2} x_{aj}^{r+1} x_{ak}^r. 				
\end{equation}
It is shown in \cite{Ni} that this data satisfies the hypotheses of the diamond lemma (although Nichols uses a slightly but not essentially different semigroup partial order), and hence the irreducible words (i.e. those which contain no subwords as in the left-hand-sides of the reductions above) form a $k$-basis for $H(M_n(k)^*)$ (with $r\in\mb N$). In \cite{Sc} it is claimed that the same is true for $H_\infty(M_n(k)^*)$ with $r\in\mb Z$, and we claim here that this holds for $H_d$ and $r\in\mb Z/2d$ as well, and hence that the irreducible words form a basis for $\tilde H$. 

As mentioned before, there are some problems when $d=1$: more ambiguities appear in this situation, which would not appear otherwise. This means that we need to check that these ambiguities are resolvable (using the language of \cite{Be}). One obvious example of such an ambiguity is $x_{nn}^r x_{nn}^{r+1}$, since for $d=1$ both \equref{red1} and \equref{red2} can be used to reduce this word. Similarly, another example (and the most tedious to resolve) of ambiguity which doesn't arise in general is $x_{in}^r x_{nn-1}^{r+1} x_{n-1 n-1}^r x_{n-1 j}^{r+1}$. Although, luckily, the new ambiguities do resolve under the reductions above, we do not perform the long but entirely straightforward calculations here. 

By a slight abuse of notation, we denote the images of the $x_{ij}^r$'s in the quotient $\tilde H$ of $k\la \mc H\ra$ by the same symbols. As noted above, the Hopf algebra structure is given by the fact that for every $r$, the $x_{ij}^r$ behave as the usual generators of an $n\times n$ matrix coalgebra (as in \equref{matco}), and the antipode acts by
\[
S(x_{ij}^r)=x_{ji}^{r+1},\ \forall r,i,j. 
\]
The coalgebra $M_n(k)^*\subset \tilde H(M_n(k)^*)$ is identified with $\{x_{ij}^0\}_{i,j}$, and hence $\{x_{ij}^r\}_{i,j}$ are its iterates through the antipode. Again, we refer to \cite{Ni, Sc} for details. 

In order to deal effectively with monomials in the $x_{ij}^r$'s, we use the following notation: bold letters such as ${\bf r}$ and ${\bf i}$ represent vectors of indices, i.e. ${\bf r}=(r_1,r_2,\ldots,r_t)$. By $x_{{\bf i}{\bf j}}^{{\bf r}}$ we mean the monomial $x_{i_1j_1}^{r_1}\ldots x_{i_tj_t}^{r_t}$. Note that the length of a vector ${\bf i}$ may vary, but in order for $x_{{\bf i}{\bf j}}^{{\bf r}}$ to make sense, ${\bf i}$, ${\bf j}$ and ${\bf r}$ must all have the same length. 

Notice also that given a fixed vector ${\bf r}$ as above, the linear span of the monomials $x_{{\bf i}{\bf j}}^{{\bf r}}$ is a subcoalgebra $C_{\bf r}$ of $\tilde H(M_n(k)^*)$. Moreover, $\tilde H(M_n(k)^*)$ is the sum of all $C_{\bf r}$'s, so every {\it simple} subcoalgebra is contained in one of the $C_{\bf r}$'s.  

The notation introduced above will be used freely throughout the rest of the paper. Here are a few more observations on free Hopf algebras which will be useful in the sequel:

\begin{remark}\relabel{scalarext}

It is shown in \cite{Ta} that the functor $H(-)$ behaves well with respect to scalar extension to a larger field. More precisely, if $k\to K$ is a field extension, then $H(C)\otimes K$ is naturally isomorphic to $H(C\otimes K)$. The analogous results for $\tilde H$ are very easy to prove, using the universal property of $\tilde H(C)$.

\end{remark}

\begin{remark}\relabel{corad}

Also in \cite{Ta}, it is shown that for any coalgebra $C$, if $C_0$ denotes its coradical and $C=C_0\oplus V$ for some vector space $V$, then $H(C)$ is $H(C_0)\coprod T(V)$ as an algebra. Here, $T(V)$ is the tensor algebra on $V$ and the coproduct is in the category of algebras. The same is in fact true if we replace $H$ with $\tilde H$. This can be seen by examining Takeuchi's proofs (\cite[Lemmas 26, 27, 28]{Ta}) and checking that they work in general.

In particular, it follows easily from this that an inclusion of coalgebras $C\to D$ induces an inclusion $\tilde H(C)\to\tilde H(D)$.  

\end{remark}

\section{Main results}\selabel{2}

As outlined in the introduction, the purpose of this section is to find all small subcoalgebras and all endomorphisms of $\tilde H(M_n(k)^*)$. The latter will be a consequence of the former, since, by the universality property of $\tilde H$, an endomorphism of $\tilde H(C)$ is the same as a coalgebra map from $C$ to $\tilde H(C)$. We will first state the main result; as before, $\tilde H(-)$ is one of $H(-),\ H_\infty(-)$ or $H_d(-)$ for some positive integer $d$.   

We have mentioned before that the cases $d=1,n=2$, positive characteristic pose some problems. It will be convenient to have a short phrase which refers to all other cases; hence, we say that we are in a {\it tame} case or situation if either $d\ge 2$, or $n>2$, or ${\rm char}(k)=0$. Otherwise, we say that we are in a {\it wild} case.

\begin{theorem}\thlabel{main}

Let $n>1$ be a positive integer. In a tame situation, the only subcoalgebras of $\tilde H(M_n(k)^*)$ of dimension $\le n^2$ are $k$ and the iterates of $M_n(k)^*\subset\tilde H(M_n(k)^*)$ by the antipode. 

\end{theorem}

\begin{remark}\relabel{chars}

In fact, it turns out that everything works fine as long as the characteristic is not $2$ or $3$. However, I've chosen to state the theorem as above, in order to keep the proof shorter (it is long enough as it is), and because it did not seem worthwhile to insist on the greatest possible generality. 

\end{remark}

Before going into the proof of the theorem, we record the desired consequences, namely the determination of the endomorphisms of $\tilde H(M_n(k)^*)$ in a tame case. By the functoriality of $\tilde H$, an (anti)endomorphism of $C$ induces an (anti)endomorphism of $\tilde H(C)$. We identify these, to avoid having to repeat the words ``induced by'' all the time.

\begin{proposition}\prlabel{end}

Let $n>1$ be a positive integer. In a tame case, the endomorphisms of $\tilde H(M_n(k)^*)$ are of one of the following types:

\begin{enumerate}
\renewcommand{\labelenumi}{(\alph{enumi})}

\item The trivial endomorphism, induced by $M_n(k)^*\stackrel{\ve}{\rightarrow} k\to\tilde H$;

\item $S^{2t}\circ\alpha$, where $t\ge 0$ and $\alpha$ is an automorphism of $M_n(k)^*$;

\item $S^{2t+1}\circ T\circ\alpha$, where $t\ge 0,\ \alpha$ is an automorphism of $M_n(k)^*$, and $T$ is the transposition map $x_{ij}^0\mapsto x_{ji}^0$ on $M_n(k)^*\subset\tilde H(M_n(k)^*)$.

\end{enumerate}

\end{proposition}

\begin{remark}\relabel{ref}

By the Skolem-Noether Theorem, the automorphisms of $M_n(k)^*$ are precisely the conjugations by $GL_n(k)$. But it is easily seen that in general, for a coalgebra $C$, a map $C\to k$ is convolution-invertible if and only if it factors through some algebra map $\tilde H(C)\to k$ (for $\tilde H(-)=H(-)$, for example, this follows immediately from \cite[$\S$2, Proposition 4]{Ta}, which characterizes algebra maps out of $H(C)$ in terms of maps out of $C$). This means that the automorphisms $\alpha$ in the statement of \prref{end} are precisely the {\it inner} automorphisms of $\tilde H(M_n(k)^*)$, in the sense that they are the conjugations (under convolution) by the algebra maps $\tilde H(M_n(k)^*)\to k$.  

\end{remark}

\renewcommand{\proofname}{Proof of \prref{end}}
\begin{proof}

As observed before, an endomorphism of $\tilde H(C)$ is determined uniquely by a coalgebra map $C\to\tilde H(C)$, so we focus on finding these. Of course, the image of a coalgebra map $C=M_n(k)^*\to \tilde H(C)$ is a subcoalgebra of dimension no larger than $n^2$, so \thref{main} applies. We thus find that our maps go either to $k$ (in which case it can only be the counit of $C$, and we are in situation (a)), or to some iterate $S^r(C)$. 

Up to an automorphism of $C=\{x_{ij}^0\ |\ i,j\}$, the map in question is $x_{ij}^0\mapsto x_{ij}^r$. This is exactly $S^r$ if $r$ is even, and $S^r\circ T$ if $r$ is odd.  
\end{proof}
\renewcommand{\proofname}{Proof}

For the proof of \thref{main}, we'll need some auxiliary results. The following lemma is an elementary linear algebra fact, whose proof we leave to the reader:

\begin{lemma}\lelabel{linalg}

Let $V,W$ be vector spaces, and $X\le V$, $Y\le W$ vector subspaces. Suppose we have an element
\[
\sum_{i=1}^p a_i\otimes b_i \in X\otimes Y,
\]
where $a_i\in V$ are linearly independent, and similarly, $b_i\in W$ are linearly independent. Then, both $X$ and $Y$ have dimension $\ge p$. 

\end{lemma}

We now prove a result in some sense weaker than the statement of \thref{main}, but which holds in wild cases too.

\begin{proposition}\prlabel{partial}

Let $n>1$ be a positive integer. The subcoalgebras of $\tilde H(M_n(k)^*)$ different from $k$ have dimension $\ge n^2$. 

\end{proposition}

\begin{proof}

By \reref{scalarext}, it suffices to consider the case when our base field $k$ is algebraically closed. This is to ensure that simple coalgebras are actually matrix coalgebras, which will be useful in the proof. Hence, throughout the rest of the argument, $k$ is assumed to be algebraically closed. 

Let $H=\tilde H(M_n(k)^*)$, and consider an arbitrary element $x\in H$. $x$ can be written as a linear combination of irreducible monomials in the standard algebra generators $x_{ij}^r$ introduced before. Let $x_{{\bf i}{\bf j}}^{{\bf r}}$ be such a monomial, having maximal length $\ell$ among the monomials appearing in $x$. We assume $\ell\ge 2$.

We look at $\Delta(x)$, using the matrix comultiplication 
\[
\Delta(x_{ij}^r)=\sum_{a=1}^n x_{ia}^r\otimes x_{aj}^r,
\]
and expanding. To get the final result in reduced monomials, we might have to reduce some of the monomials we get by expansion. 

Now fix a pair $(\alpha,\beta)$ of distinct indices in $\ol{1,n}$, and consider the vector ${\bf u}=(\alpha,\beta,\alpha,\beta,\ldots)$, having the same length $\ell$ as ${\bf r}$. The term $x_{{\bf i}{\bf u}}^{{\bf r}}\otimes x_{{\bf u}{\bf j}}^{{\bf r}}$ will appear in $\Delta(x)$ (after all the reductions have been made). This follows because on the one hand $x_{{\bf i}{\bf u}}^{{\bf r}}$ and $x_{{\bf u}{\bf j}}^{{\bf r}}$ are reduced (as a consequence of the fact that $x_{{\bf i}{\bf j}}^{{\bf r}}$ was reduced and the form of the reduction rules \equref{red1}, \equref{red2}, etc.), and on the other hand because of the maximality of the length of $x_{{\bf i}{\bf j}}^{{\bf r}}$, which implies that $x_{{\bf i}{\bf u}}^{{\bf r}}$ and $x_{{\bf u}{\bf j}}^{{\bf r}}$ cannot be non-trivial reductions of some other monomials we run into when trying to compute $\Delta(x)$. 
 
Applying this argument to all $n(n-1)$ ordered pairs $(\alpha,\beta)$ of distinct indices in $\ol{1,n}$, we find that $\Delta(x)$, in its unique form as a linear combination of tensor products of reduced monomials, contains all $x_{{\bf i}{\bf u}}^{{\bf r}}\otimes x_{{\bf u}{\bf j}}^{{\bf r}}$ (as before, having fixed the pair $(\alpha,\beta)$, we set ${\bf u}$ to be $(\alpha,\beta,\alpha,\ldots)$). 
 
Now consider a simple subcoalgebra $C\subset H$. $C$ is a matrix coalgebra, because $k$ is algebraically closed. Assuming $C$ is neither $k$ nor one of the coalgebras $\{x_{ij}^r\}_{i,j}$, all of its elements contain monomials of length $\ge 2$. In conclusion, the argument above applies to all elements $x\in C$. If $C$ is, say, an $m\times m$ matrix coalgebra, then we can find an element $x\in C$ such that $\Delta(x)\in M\otimes N$, where $M,N$ are linear spaces of dimension $m$ (for example, $x$ can be part of a system of matrix generators, such as the $t_{ij}$ in \equref{matco}). Pick a monomial $x_{{\bf i}{\bf j}}^{{\bf r}}$ for $x$, using the notations above. Let $P\le H$ (resp. $Q\le H$) be the linear subspace generated by all monomials not of the form $x_{{\bf i}{\bf u}}^{{\bf r}}$ (resp. $x_{{\bf i}{\bf j}}^{{\bf r}}$), where ${\bf u}$, as before, ranges through $(\alpha,\beta,\alpha,\ldots)$. Now, applying \leref{linalg} to $V=H/P$, $W=H/Q$, $X=M+P/P$ and $Y=N+Q/Q$, we conclude that the dimension of $M+P/P$ (and hence that of $M$) is at least $n(n-1)$. Hence, $m\ge n(n-1)\ge n$. 

This proves that simple subcoalgebras of $H$ are either $k$, one of the iterates of $\{x_{ij}^0\}$ through the antipode, or $m\times m$ matrix subcoalgebras with $m\ge n(n-1)\ge n$. This implies that any subcoalgebra $C$ of $H$ of dimension $\le n^2$ is either connected with coradical $k$, or an $n\times n$ matrix coalgebra. The former is impossible, however, unless $C=k$, because $H$ has no non-zero primitive elements: such a primitive element $x$ would have to contain a monomial of length $\ge 2$, and the argument above would show that the image of $\Delta(x)$ in $(H/P)\otimes(H/Q)$ is non-zero; but $1\in P\cap Q$, so the image of $1\otimes x+x\otimes 1$ is $(H/P)\otimes(H/Q)$ is zero.

This finishes the proof of the proposition. 
\end{proof}

The argument used in the previous proof is the essential ingredient in \thref{main}, and it will appear in various guises throughout the rest of our proof of the main theorem. As a consequence of this argument, we already have the following partial result:

\begin{corollary}\colabel{main n>2}

The conclusion of \thref{main} holds if $n>2$. 

\end{corollary}

\begin{proof}

Again, we may as well assume the base field $k$ is algebraically closed. 

We remarked in the last paragraph of the proof for \prref{partial} that (a) the simple subcoalgebras of $\tilde H(M_n(k)^*)$ different from $k$ or $\{x_{ij}^r\}_{i,j}$ are $m\times m$ matrix coalgebras with $m\ge n(n-1)$ (which is {\it strictly} larger than $n$ if $n>2$) and (b) there are no subcoalgebras with coradical $k$. The conclusion is now clear. 
\end{proof}

In view of this corollary, we can focus on the case $n=2$, although the simplification is only notational. The following lemma will also come in handy:

\begin{lemma}\lelabel{2x2}

Suppose $n=2$, and let ${\bf r}=(r_1,r_2,\ldots)$ be a vector of length at least $2$, of elements of $\mb N$, $\mb Z$ or $\mb Z/2d$ according as $\tilde H(-)$ is $H(-)$, $H_\infty(-)$ or $H_d(-)$, respectively. In a tame case, the linear span $D=D_{\bf r}$ of the four elements $x_{{\bf i}{\bf j}}^{{\bf r}}$ where ${\bf i}$ and ${\bf j}$ are alternating vectors of the form $(1,2,1,\ldots)$ or $(2,1,2,\ldots)$ is not a subcoalgebra of $\tilde H(M_2(k)^*)$.   

\end{lemma}

\begin{proof}

Remember that by our conventions at the beginning of this section, and considering that $n=2$, being in a tame case means that either $\tilde H(-)$ is not $H_1(-)$, or that we are working in characteristic zero. 

Assume first that $\tilde H(-)$ is not $H_1(-)$. This means that we can find consecutive entries $r_i$ and $r_{i+1}$ of the vector ${\bf r}$ (which has length at least $2$, by the hypothesis) such that either $r_{i+1}\ne r_i+1$, or $r_{i+1}\ne r_i-1$. To fix ideas, suppose, for example, that $r_2\ne r_1+1$; the general case is entirely analogous.

Just as before, in the proof of \prref{partial}, we are going to try to compute $\Delta(x_{{\bf i}{\bf j}}^{{\bf r}})$ for some monomial in $D$ by using the matrix comultiplication rules and expanding. Let ${\bf u}$ be the vector $(2,2,1,2,\ldots)$ of the same length as ${\bf r}$. ${\bf u}$ has $2$ as its first entry, and then alternates, starting with $2$ again. Because of our assumption on ${\bf r}$, the monomial $x_{{\bf i}{\bf u}}^{{\bf r}}$ is reduced (this is easily seen by examining the reduction rules \equref{red1}-\equref{red4}). Moreover, the same reduction rules imply that $x_{{\bf i}{\bf u}}^{{\bf r}}$ cannot be obtained as a non-trivial reduction from another monomial appearing in our computation of $\Delta(x_{{\bf i}{\bf j}}^{{\bf r}})$. It follows then that after reducing everything in the expression of $\Delta(x_{{\bf i}{\bf j}}^{{\bf r}})$, there will be at least one term of the form $\pm x_{{\bf i}{\bf u}}^{{\bf r}}\otimes\bullet$ left. This term is not an element of $D\otimes D$ (because ${\bf u}$ is not alternating), and we are done. 

Now assume $\tilde H(-)=H_1(-)$, but ${\rm char}(k)=0$. Because the entries of the vector ${\bf r}$ are elements of $\mb Z/2$, there's no difference now between $r_i+1$ and $r_i-1$. The previous argument still works if two consecutive entries of ${\mb r}$ are equal, but not if ${\mb r}$ is an alternating vector (i.e. any two consecutive entries are different). Nevertheless, we try to apply the same technique, and compute $\Delta(x_{{\bf i}{\bf j}}^{{\bf r}})$ for some reduced monomial in $D$. 

Let ${\bf u}$ be the vector $(1,1,1,\ldots)$, of the same length as ${\bf r}$. Notice that by the reduction rules \equref{red1} - \equref{red4}, for any vector ${\bf v}$ of the same length, the coefficient of $x_{{\bf i}{\bf u}}^{{\bf r}}$ in the reduced form of $x_{{\bf i}{\bf v}}^{{\bf r}}$ is equal to the coefficient of $x_{{\bf u}{\bf j}}^{{\bf r}}$ in the reduced form of $x_{{\bf v}{\bf j}}^{{\bf r}}$. It follows then, because we are working in characteristic zero, that after performing all the reductions, the coefficient of $x_{{\bf i}{\bf u}}^{{\bf r}}\otimes x_{{\bf u}{\bf j}}^{{\bf r}}$ in $\Delta(x_{{\bf i}{\bf j}}^{{\bf r}})$ is positive. In particular, $\Delta(x_{{\bf i}{\bf j}}^{{\bf r}})$ does not belong to $D\otimes D$. 
\end{proof}

\begin{remark}\relabel{sec4}

We will see below, in \seref{4}, that the tame case hypothesis is necessary. 

\end{remark}

Finally, we are ready now to finish the proof of the main theorem.

\renewcommand{\proofname}{Proof of \thref{main}}
\begin{proof}

As remarked repeatedly before, we can assume the base field is algebraically closed. We already know, from the proof of \prref{partial}, that (for the purpose of our theorem) it suffices, over an algebraically closed field, to look only at matrix subcoalgebras of $H=\tilde H(M_n(k)^*)$. Also, we assume $n=2$, as permitted by \coref{main n>2}. Finally, by an observation made at the end of \seref{1}, a matrix subcoalgebra of $H$ is contained in some $C_{\bf r}$, the linear span of the monomials $x_{{\bf i}{\bf j}}^{\bf r}$ for some fixed ${\bf r}$. 

In line with the previous paragraph, let $C\subseteq C_{\bf r}\subset H$ be an $m\times m$ matrix subcoalgebra of $H$, with $m\le n$. We may as well assume that the length $\ell$ of ${\bf r}$ is at least $2$. Pick an $x\in C$, and let $x_{{\bf i}{\bf j}}^{{\bf r}}$ be a reduced monomial appearing in $x$. We saw in the proof for \prref{partial} that after performing all the reductions, $\Delta(x_{{\bf i}{\bf j}}^{{\bf r}})$ contains both terms of the form $x_{{\bf i}{\bf u}}^{{\bf r}}\otimes x_{{\bf u}{\bf j}}^{{\bf r}}$, where ${\bf u}$ is one of the two alternating vectors of length $\ell$ (either $(1,2,1,\ldots)$ or $(2,1,2,\ldots)$).

The proof for \prref{partial} (more specifically the part of the proof which used \leref{linalg}, contained in the last two paragraphs of the proof) also shows that if $\Delta(x_{{\bf i}{\bf j}}^{{\bf r}})$ were to contain $x_{{\bf a}{\bf b}}^{{\bf s}}\otimes x_{{\bf c}{\bf d}}^{{\bf t}}$ with neither ${\bf b}$ nor ${\bf c}$ alternating of length $\ell$, then $x$ could not be one of the matrix generators of $C$. It follows that for any such generator, all the terms of $\Delta(x)$ (after all the reductions have been made) are multiples either of $x_{{\bf i}{\bf u}}^{{\bf r}}\otimes\bullet$ or $\bullet\otimes x_{{\bf u}{\bf j}}^{{\bf r}}$, with ${\bf u}$ alternating of length $\ell$. But because of the maximality of the length of $x_{{\bf i}{\bf j}}^{{\bf r}}$ in $x$, it's clear that the only possible such terms are the multiples of $x_{{\bf i}{\bf u}}^{{\bf r}}\otimes x_{{\bf u}{\bf j}}^{{\bf r}}$ (in other words, if $x_{{\bf i}{\bf j}}^{{\bf r}}\otimes\bullet$ were to appear in $\Delta(x)$, the only possibility for $\bullet$ would be $x_{{\bf u}{\bf j}}^{{\bf r}}$). 

In conclusion, for matrix generators $x$ of $C$, $\Delta(x)$ is a linear combination of the two $x_{{\bf i}{\bf u}}^{{\bf r}}\otimes x_{{\bf u}{\bf j}}^{{\bf r}}$, with ${\bf u}$ alternating of length $\ell$. But by using the counit identities on $x$, we see that this implies that $x$ is a member of what in the statement of \leref{2x2} was denoted by $D_{\bf r}$, and hence that our $C\subseteq C_{\bf r}$ be $D_{\bf r}$. But \leref{2x2} says precisely that in a tame case, $D_{\bf r}$ is not a subcoalgebra (for ${\bf r}$ of length $\ge 2$). This finishes the proof of the theorem.
\end{proof}
\renewcommand{\proofname}{Proof}

As a final remark, we record the following consequence of \thref{main}:

\begin{corollary}\colabel{main}

Let $n>1$ be a positive integer. In a tame case, the only right comodules over $H=\tilde H(M_n(k)^*)$ of dimension $\le n$ are (a) the direct sums of $\le n$ copies of the trivial comodule, and (b) the iterated duals of the $n$-dimensional comodule obtained by scalar corestriction from $M_n(k)^*\to H$.  

\end{corollary}

\begin{proof}

Let $M$ be a right comodule over $H$, of dimension $m\le n$, with comodule structure map $\rho:M\to M\otimes H$. If $e_i,\ i=\ol{1,m}$ is a basis for $M$, then we get elements $c_{ij}$ of $H$ by
\[
\rho_{e_j}=\sum_i e_i\otimes c_{ij}.  
\]
It's easy to see that the $c_{ij}$ satisfy matrix coalgebra-type relations, as in \equref{matco}, or, in other words, we have a coalgebra map from $M_n(k)^*$ to $H$ sending the standard generators $t_{ij}$ to $c_{ij}$. But this means that the $c_{ij}$ form a subcoalgebra of $H$ of dimension $\le n^2$, and the conclusion follows immediately from \thref{main}. 
\end{proof}

\section{Centers of some categories}\selabel{3}

Here, as an application of \thref{main}, we determine the centers of the categories $\halg$, $\halg_\infty$ and $\halg_d$. Because these centers are all monoids with a ``multiplicative zero'', namely the natural transformation which is given on each Hopf algebra (or Hopf algebra with bijective antipode, or Hopf algebra $H$ with $S_H^{2d}=\{\rm id\}$) by the composition between the unit and the counit, it will be convenient to have a notation for this phenomenon. Hence, we introduce the following notation:

For a monoid $M$, denote by $M^+$ the monoid which as a set is $M\cup\{0\}$, with multiplication defined by the one in $M$ and by 
\[
0x=x0=0,\ \forall x\in M. 
\]
In other words, $M^+$ is obtained from $M$ by appending a multiplicative zero. 

In the following statement, $\mb N,\mb Z$ and $\mb Z/2d$ are monoids with their usual additive structure. Notice that in each of our categories, there is an endo-natural transformation of the identity functor given by the square of the antipode on each object of the category. To avoid cumbersome language, we refer to this natural transformation as {\it being} the square of the antipode.

\begin{theorem}\thlabel{centers}

The centers of $\halg$, $\halg_\infty$ and $\halg_d$ are $\mb N^+$, $\mb Z^+$, and $(\mb Z/2d)^+$, respectively, where $\mb N$, $\mb Z$ and $\mb Z/2d$ are generated by the square of the antipode. In all three cases, the multiplicative zero is given by the trivial endomorphism.     

\end{theorem}

\begin{proof}

We prove the statement for $\halg$; the proofs in the other two cases are entirely parallel. 

First, looking at the action of the antipode on our elements $x_{ij}^0$ in some $H(M_n(k)^*)$, it's clear that the different powers of the antipode induce different endo-natural transformations of the identity, and hence the monoid generated by $S^2$ and the trivial endomorphism is indeed $\mb N^+$. The interesting part is showing that conversely, {\it every} element of the center is either trivial or induced by some even power of the antipode.  

Let $\eta$ be an endo-natural transformation of the identity functor on $\halg$. This means that for every Hopf algebra $H$, we are given an endomorphism $\eta_H$ of $H$ such that 
\begin{diagram}
H            &\rTo^{\eta_H}        &H       \\
\dTo<f       &                     &\dTo>f  \\
K            &\rTo^{\eta_K}        &K
\end{diagram}
commutes for every Hopf algebra map $f:H\to K$. Let us look at what $\eta_H$ might be for $H=H(M_n(k)^*)$ for some fixed $n>1$ (we would take $n>2$ if we were dealing with $\halg_1$ instead of $\halg$, to make sure we are in a tame situation). \thref{main} says that there are three cases: 

(1) $\eta_H$ is of the form $S^{2t}\circ\alpha$ for some automorphism $\alpha$ of $M_n(k)^*$. It is clear (for example from the structure of the basis of $H$ we've been working with) that the map ${\rm Aut}(M_n(k)^*)\to {\rm End}(H(M_n(k)^*))$ given by $\beta\mapsto S^{2t}\circ\beta$ is injective. From this and the commutativity of the diagram above for $\eta_H=S^{2t}\circ \alpha$ and $f=\beta\in{\rm Aut}(M_n(k)^*)$ it follows that $\alpha$ is in the center of ${\rm Aut}(M_n(k)^*)$. This implies $\alpha={\rm id}$, and hence $\eta_H$ is an even power of the antipode. 

(2) $\eta_H$ is of the form $S^{2t+1}\circ T\circ \alpha$, where $T$ is the transposition on $M_n(k)^*$, and $\alpha$ is an automorphism of the matrix coalgebra. Just as before, consider our commutative diagram with $K=H$ and $f=\beta$, some automorphism of $M_n(k)^*$. The same argument as in (1) (and the fact that every endomorphism of $H$ commutes with the antipode) shows that $T\alpha\beta=\beta T\alpha$ for arbitrary $\beta$. This is easily seen to be impossible for $n>2$, and hence (2) is ruled out. 

(3) $\eta_H$ is the trivial endomorphism of $H$. 

Denote $\eta_H,\ H=H(M_n(k)^*)$ by $\eta_n$. We now know that $\eta_n$ is either $S^{2r}$ for some $r$ or trivial (for $n>2$, at least). I claim that either we have the same $r$ for all $n$, or $\eta_n$ is trivial for all $n$. First, notice that the claim finishes the proof. To see this, suppose, for example, that $\eta_n=S^{2r}$ for every large $n$ (the case where $\eta_n$ are all trivial is similar). Now, by the commutativity of the square diagram above, $\eta_K$ is going to be $S_K^{2r}$ for every quotient of a Hopf algebra of the form $H(M_n(k)^*)$ for large $n$. But on the one hand, every finite-dimensional coalgebra is a quotient of some $M_n(k)^*$, and on the other hand, every Hopf algebra is the union of its finite-dimensional subcoalgebras; this implies that every Hopf algebra is a union of quotients of Hopf algebras of the form $H(M_n(k)^*)$, and we are done. 

All that remains is to prove the claim. Say for some fixed $n>2$, $\eta_n$ is $S^{2r}$, while $\eta_{n+1}$ is $S^{2s}$ (again, the case when one of $\eta_n$, $\eta_{n+1}$ is trivial is analogous). This means, in terms of our standard algebra generators $\{x_{ij}^r\}$ for $H(M_n(k)^*)$ and $\{y_{ij}^r\}$ for $H(M_{n+1}(k)^*)$, that $\eta_n$ is the endomorphism induced by $x_{ij}^0\mapsto x_{ij}^{2r}$, while $\eta_{n+1}$ is induced by $y_{ij}^0\mapsto y_{ij}^{2s}$.  

Now let $C$ be the quotient of $M_{n+1}(k)^*$ by the coideal spanned by $y_{n+1j}^0,\ j=\ol{1,n}$. We denote the images of $y_{ij}^0$ in $C$ by the same symbols. $\eta_H$ will be $S^{2s}$ for $H=H(C)$. At the same time, however, we have an inclusion $H(M_n(k)^*)\to H(C)$ (\reref{corad}) given by the inclusion $M_n(k)^*\to C$ given as $x_{ij}^0\mapsto y_{ij}^0,\ i,j=\ol{1,n}$. It follows now that $\eta_n$ is both $S^{2r}$ and $S^{2s}$. As the $S^i$ are different for different $i$ on $H(M_n(k)^*)$ (by looking at how the powers of the antipode act on the $x_{ij}^0$), we get $r=s$, as desired. 
\end{proof}

\section{What about $H_1(M_2(k)^*)$ in positive characteristic?}\selabel{4}

The purpose of this short section is to point out that, as mentioned several times before, the tame case hypothesis in \thref{main} is actually necessary. More specifically, we have counterexamples in characteristics $2$ and $3$. We observed in \reref{chars} that in fact \thref{main} works even for $H_1(M_2(k)^*)$ in positive characteristic as long as it is different from $2$ or $3$, but we will not prove this here. The proof consists of making a slightly more detailed analysis of what can go wrong with the arguments in \seref{2}, using the same techniques as before.

\begin{example}\exlabel{char2}

Suppose the base field $k$ has characteristic $2$, and let ${\bf r}$ be either $(0,1)$ or $(1,0)$, where $0,1$ are the elements of $\mb Z/2$. Then, using the notation from \leref{2x2}, $D_{\bf r}$ is a $2\times 2$ matrix subcoalgebra of $H_1(M_2(k)^*)$. 

\end{example}

\begin{proof}

This is a simple verification. Assume for example that ${\bf r}$ is $(0,1)$. We check that $\Delta(x_{11}^0x_{22}^1)$ does indeed belong to $D_{\bf r}\otimes D_{\bf r}$, and leave the rest to the reader. 

We have 
\begin{align}\eqlabel{char2}
\Delta(x_{11}^0x_{22}^1)    &=\Delta(x_{11}^0)\Delta(x_{22}^1) = (x_{11}^0\otimes x_{11}^0+x_{12}^0\otimes x_{21}^0)
(x_{21}^1\otimes x_{12}^1+x_{22}^1\otimes x_{22}^1)                               \nonumber                         \\
                            &=x_{11}^0x_{21}^1\otimes x_{11}^0x_{12}^1+x_{11}^0x_{22}^1\otimes x_{11}^0x_{22}^1
                             +x_{12}^0x_{21}^1\otimes x_{21}^0x_{12}^1+x_{12}^0x_{22}^1\otimes x_{21}^0x_{22}^1.                                        
\end{align}
Now simply notice that because of the two reduction rules \equref{red1} and \equref{red2}, we have (regardless of the characteristic) 
\begin{align*}
x_{12}^0x_{22}^1       &=-x_{11}^0x_{21}^1\\
x_{21}^0x_{22}^1       &=-x_{11}^0x_{12}^1.  
\end{align*}
Because ${\rm char}(k)=2$, the first and last term in \equref{char2} cancel out. 
\end{proof}

Similarly, we have

\begin{example}\exlabel{char3}

Suppose ${\rm char}(k)=3$, and ${\bf r}$ is one of the alternating vectors $(0,1,0)$ or $(1,0,1)$ with entries from $\mb Z/2$. Then, $D_{\bf r}$ is an $2\times 2$ matrix subcoalgebra of $H_1(M_2(k)^*)$. 

\end{example}

\section*{Acknowledgement}

The author wishes to thank the referee, to whom the observatnion in \reref{ref} is due, for this and other suggestions on how to improve the manuscript.  




\begin{thebibliography}{99}



\bibitem[A]{A}
Abe, E. - Hopf algebras, Cambridge University Press 1980


\bibitem[Be]{Be}
Bergman, G. - The diamond lemma for ring theory, {\it Adv. Math.} {\bf 29} (1978), pp. 178 - 218


\bibitem[Mo]{Mo}
Montgomery, S. - Hopf algebras and their actions on rings, vol. 82 of {\it CBMS Regional Conference Series in Mathematics}, AMS, Providence, Rhode Island 1993


\bibitem[Ni]{Ni}
Nichols, W. D. - Quotients of Hopf algebras, {\it Comm. Algebra} {\bf 6} (1978), pp. 1789 - 1800


\bibitem[Sc]{Sc}
Schauenburg, P. - Faithful flatness over Hopf subalgebras: Counterexamples, appeared in {\it Interactions between ring theory and representations of algebras: proceedings of the conference held in Murcia, Spain}, CRC Press (2000), pp. 331 - 344


\bibitem[Sw]{Sw}
Sweedler, M. E. - Hopf algebras, Benjamin New York 1969 


\bibitem[Ta]{Ta}
Takeuchi, M. - Free Hopf algebras generated by coalgebras, {\it J. Math. Soc. Japan} {\bf 23} (1971), pp. 561 - 582



\end{thebibliography}
\end{document}